\documentclass{article}[12pt]
\usepackage{amsmath}
\usepackage{amssymb}
\usepackage{amscd}
\usepackage{amsthm}
\usepackage{color}
\usepackage{enumerate}
\usepackage{caption}
\usepackage{subcaption}
\usepackage{graphics,epsfig}
\usepackage[top=1in, bottom=1in, left=1in, right=1in]{geometry}
\usepackage{graphicx}
\usepackage{authblk}
\usepackage{afterpage}

\begin{document}

\title{(Un)predictability of strong El Ni\~no events}
\author{John Guckenheimer}
\affil{Mathematics Department, Cornell University, Ithaca, NY 14853}
\author{Axel Timmermann}
\affil{IBS Center for Climate Physics, Pusan National University, Busan, South Korea} 
\author{Henk Dijkstra}
\affil{Institute for Marine and Atmospheric Research Utrecht, 
Department of Physics and Astronomy, Utrecht University, Utrecht, 
The Netherlands} 
\author{Andrew Roberts}
\affil{Cerner Corporation, Kansas City, Missouri}

\date{\today}

\maketitle

\begin{abstract}
 
The El Ni\~no-Southern Oscillation (ENSO) is a mode of interannual variability in 
the coupled  equatorial Pacific coupled atmosphere/ocean system. El Ni\~no describes 
a state in which sea surface temperatures in the eastern Pacific increase 
and upwelling of colder, deep waters diminishes. El Ni\~no events typically 
peak in boreal winter, but their strength varies irregularly on decadal time 
scales. There were exceptionally strong El Ni\~no events in 1982-83, 
1997-98 and 2015-16 that affected weather on a global scale. Widely publicized 
forecasts in 2014 predicted that the 2015-16 event would occur a year earlier.
Predicting the strength of El Ni\~no is a matter of practical concern 
due to its effects on hydroclimate and agriculture around the world.
This paper discusses the frequency and regularity of strong El Ni\~no events
in the context of chaotic dynamical systems. We discover a  mechanism that 
limits their predictability in a conceptual ``recharge oscillator'' model of ENSO. 
Weak seasonal forcing  or noise
in this model can induce irregular switching  between an oscillatory state  
that has strong  El Ni\~no events  and a chaotic state that lacks strong events,
In this regime, the timing of strong  El Ni\~no events on decadal time scales is 
unpredictable.
\end{abstract}

 \section{Introduction}

The problem in predicting  El Ni\~no events is that they occur quite  irregularly, and their 
development seems to be different each time \cite{McPhaden2015}. In addition to observations, the amplitudes of El Ni\~no events in a long simulation of the GFDL CM2.1 coupled ocean-atmosphere General Circulation model with constant forcing are highly variable on decadal time scales \cite{Wittenberg2009}.
The predictability horizon of individual El Ni\~no events is thought to be about 6-9 months, depending on the season and the phase of the ENSO cycle.
During boreal spring the  coupled ocean-atmosphere system is thought to be at its frailest state 
\cite{Webster1992QJRMS}.  Then the system is most susceptible to 
perturbations \cite{Webster1995MAP} which leads  to a `spring' predictability 
barrier in April/May \cite{Latif1994CD}.  The role of the  initial error pattern has been 
emphasized and in particular its interaction with  the seasonal cycle and the internal ENSO
cycle  \cite{Mu2007a,Duan2009,Yu2012}. 
However, since detailed surface and ocean subsurface observations of ENSO began in the 1950s, very strong El Ni\~no events 
have occurred only once every 15-20 years (1982, 1997 and 2015)
and the processes controlling their predictability horizon are largely unknown.  

Similar to 
investigations of tipping points \cite{Scheffer2009a}, we would like to identify precursors that 
predict the strength of impending events and their frequency on both shorter and longer 
time scales. 
Processes determining ENSO variability  and its limits of predictability have been 
studied in conceptual models \cite{Suarez1988,Tziperman1998,Jin1997}. Such 
models seek to identify key dynamical processes in the complex, coupled 
ocean-atmosphere system in the tropical Pacific.  

This paper discovers a new 
dynamical mechanism that limits the predictability 
of strong  El Ni\~no events in the low-dimensional conceptual Jin-Timmermann (JT) 
model of ENSO  originally proposed  by Jin \cite{Jin1997} and then extended  by  
Timmermann et al. \cite{tim03}.  The state 
space variables of this deterministic  model are sea surface temperatures of the 
equatorial western Pacific and  eastern Pacific, and the thermocline depth (e.g. the 
depth of the 20$^\circ$C isotherm) of the 
western Pacific. Nonlinear  terms in the model  are associated with the wind-induced anomalous advection of sea surface temperature anomalies and thermocline dynamics that  affect these state variables. 
Comparisons of this reduced order model with  observational data confirm that the model embodies  the basic features of ENSO in the tropical 
Pacific \cite{an_jin2004}. 

Recently,  a  multiple time scale analysis of the JT model  was performed to 
gain new  insight  into its dynamics  \cite{roberts2016}.  This work transformed the model  
equations into an equivalent dimensionless system whose equations are
\begin{equation}  \label{fast}
\begin{cases}
x' &=\rho \delta (x^2 - ax) + x \left( x+y+c - c \tanh (x+z) \right)\\
y' &= -\rho \delta (ay + x^2) \\
z' &= \delta (k-z- \frac{x}{2}), 
\end{cases}
\end{equation}
The variables $x,y$ and $ z$ of this  model represent the sea surface 
temperature difference between the eastern and western equatorial Pacific, the sea 
surface temperature of the western equatorial Pacific relative to a nominal 
reference temperature, and the  thermocline depth of the western Pacific, respectively. 
Each of the five 
parameters $\delta,\rho,c,k$ and $a$ represents a combination of 
physical characteristics of the tropical Pacific. The tuning of these 
parameters and the sensitivity of the El Ni\~no variability  to these 
parameters is discussed in \cite{roberts2016}. Before we present additional analysis of this model, we recall relevant background material from dynamical systems theory.

\section{Sensitive Dependence to Initial Conditions and Predictability}

Strong El Ni\~no events are typically characterized by a threshold of average sea surface temperature anomaly in the eastern equatorial Pacific, e.g., the Nino 3.4 region  (5N-5S, 170W-120W).~\cite{Bamston1997} Within the setting of the JT model, we set the threshold $x > -1.5$ for strong El Ni\~no events since $x$ represents the difference in sea surface temperature of the eastern and western tropical Pacific. The events are then \emph{recurrences} to the phase space region $x> -1.5$. Predictability of strong El Ni\~no events in the JT model is then a matter of the frequency and regularity of recurrences to this region. We use a simple discrete time dynamical system to illustrate the type of unpredictability that we subsequently find in the JT model.

Sensitive Dependence to Initial Conditions \cite{Ruelle} has been a key concept in the 
study of chaotic dynamical systems: pairs of nearby initial conditions yield trajectories
that separate from each other. Statistical  theories of chaotic attractors \cite{Young} 
emphasize asymptotic properties like invariant measures, Lyapunov exponents and entropies. 
There are dynamical systems in which multiple time scales are an emergent phenomenon. Some trajectories spend long times near metastable invariant sets, occasionally making transitions between them. We use the term \emph{epoch} to denote a maximal time interval when the system is close to one of the metastable invariant sets. 

Our viewpoint is strongly influenced by the Wentzell-Freidlin theory of stochastically perturbed dynamical systems~\cite{Freidlin}. Imagine a deterministic dynamical system with several attractors and a small stochastic or random perturbation of this system. Wentzell and Freidlin introduce a Markov chain whose states are the attractors of the deterministic system with transition probabilities given by the observed frequencies of transitions between their basins of attractions in trajectories of the stochastic system. This same Markov model can be used for  deterministic systems which have an attractor comprised of several \emph{almost invariant} sets with infrequent transitions from one to another~\cite{Froyland}. 

\begin{figure}
\begin{center}
\includegraphics[width=0.9\textwidth]{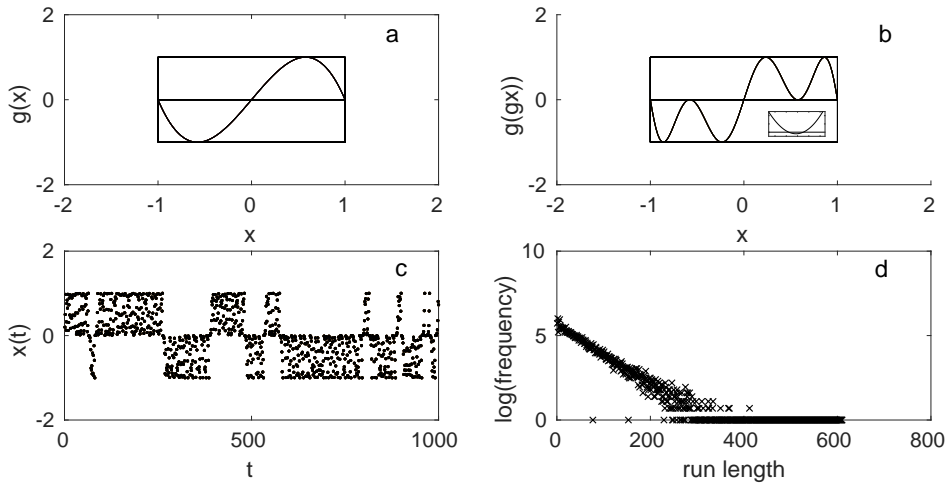}
\end{center}
  \caption{Mode switching in a discrete one dimensional iteration: (a) Graph of the function $g(x) = 2.6*x(1-x^2)$. The intervals $[g(-1/\sqrt{3}),0]$ and $[0,g(1/\sqrt{3})]$ are almost invariant, but short subintervals near $\pm 1/\sqrt{3}$ map to the opposite interval. (b) The second iterate of $g$ together with an expanded plot near $1/\sqrt{3}$ showing that the local minimum of the map is negative. (c) 1000 iterates beginning at $x=0.2$. Epochs in which the sign of the iterates remain constant are clearly visible. (d) log(number of epochs of length n) against n in a trajectory segment of $10^{6}$ iterates.}
\label{oned_mode_switch}
\end{figure}

We illustrate these concepts with iterations of the one dimensional map $g(x) = \alpha*x(1-x^2), \, x \in [-\infty,\infty]$. When $2<\alpha<3\sqrt{3}/2 \approx 2.598$, this map has two attractors in the intervals $[-1,0]$ and $[0,1]$. At $\alpha = 1/\sqrt{3}$, the two attractors merge. Iterations of one dimensional maps are an intricate subject with deep mathematical structure.~\cite{MeloStrien}. For some values of the parameter $\alpha$, $g$ will have stable periodic orbit(s) that are its attractors, but for a positive measure set of $\alpha$, the map has chaotic attractor(s). Figure~\ref{oned_mode_switch} illustrates the dynamics when $\alpha=2.6$, seemingly a value that gives rise to a single chaotic attractor. The intervals $[g(-1/\sqrt{3}),0]$ and $[0,g(1/\sqrt{3})]$ are almost invariant: trajectories spend quite long periods of time in these intervals before switching to the other interval. This is shown in Figure 1c where a trajectory segment of length $1000$ is plotted. There are epochs between mode switches of quite different lengths in this subplot. Figure 1d shows a histogram (on a log scale) of the number of epochs of different lengths in a trajectory of $10^{6}$ iterates. The longest epoch has length 610; the shortest has length 8. The approximately linear slope of the graph suggests that the epoch lengths can be modeled by a Poisson distribution. In a longer trajectory of $10^{7}$ iterates, we find that \emph{all} pairs of integers $(m,n)$ with $8 \le m,n \le 60$ occur as successive epoch lengths, indicating that the length of an epoch does not enable one to predict the length of the next epoch. Superficially,  successive epoch lengths seem to be statistically independent of each other.

The unpredictability of epoch lengths in this model is related to its sensitive dependence to initial conditions. Nearby initial conditions yield trajectories that separate from one another so that information about their approximate location is lost after a moderate number of iterates. In this example, recurrences to the intervals $[-\infty,-1]$ and $[1,\infty]$ mark the transitions from one basin to the opposite one. The set of initial conditions where a transition occurs at the $n^{th}$ iterate is a finite union of intervals, but the number of intervals increases geometrically with $n$ and, for a positive measure set of $\alpha$,  their distribution approaches an absolutely continuous invariant measure~\cite{MeloStrien}. 

We also computed the epoch lengths for sample trajectories of the stochastic map $g_r(x) = 2.595*x(1-x^2) + 0.005*\xi$ (where $\xi$ is a normally distributed random variable) to reside in one of the two attractor basins. Here, we are certain that the dynamics have an absolutely continuous invariant measure and exponential decay of correlations~\cite{Young92a}.  The resulting distribution of epoch lengths has no apparent differences from those shown in Figure~\ref{oned_mode_switch}(d) for the deterministic system with two almost invariant sets. Both the deterministic and random versions of this model have epoch lengths that are unpredictable from past history.

\section{Two Attractors of the JT Model}

\begin{figure}
\includegraphics[width=0.45\textwidth]{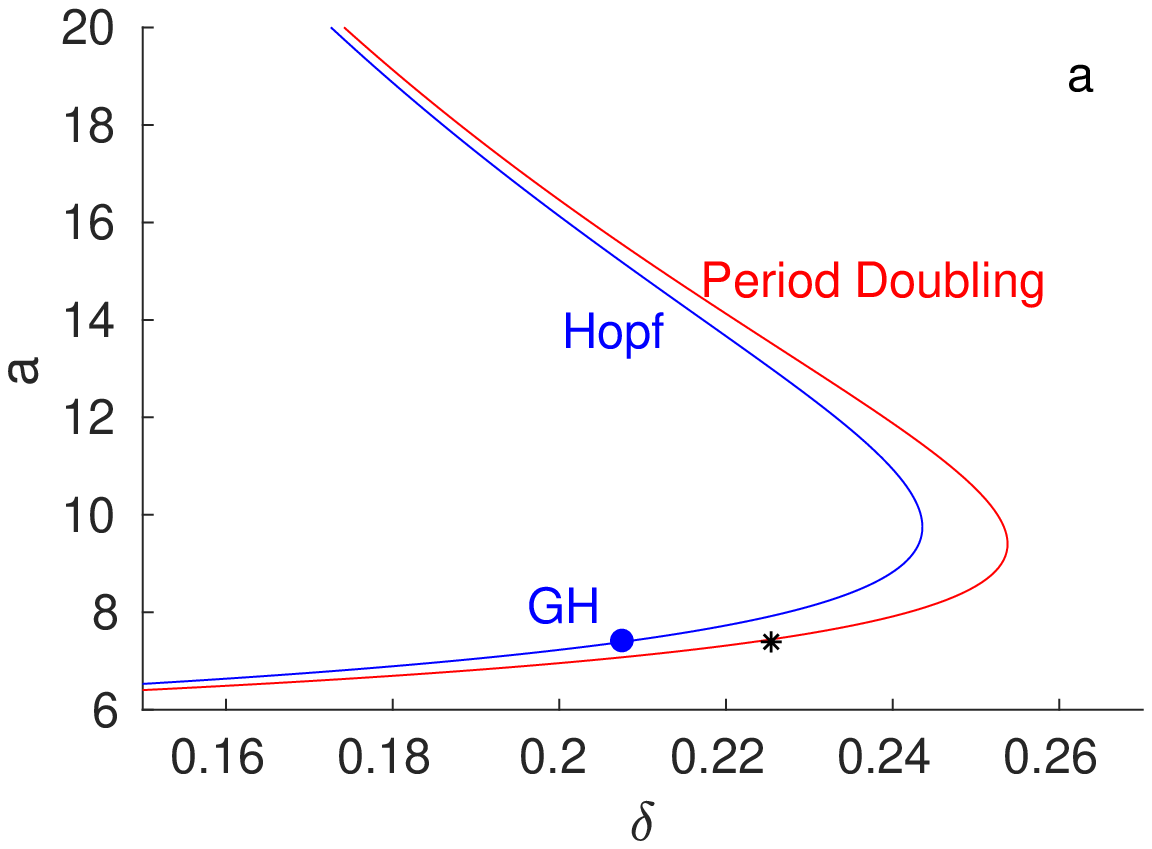}
\hspace{0.05\textwidth}
\includegraphics[width=0.45\textwidth]{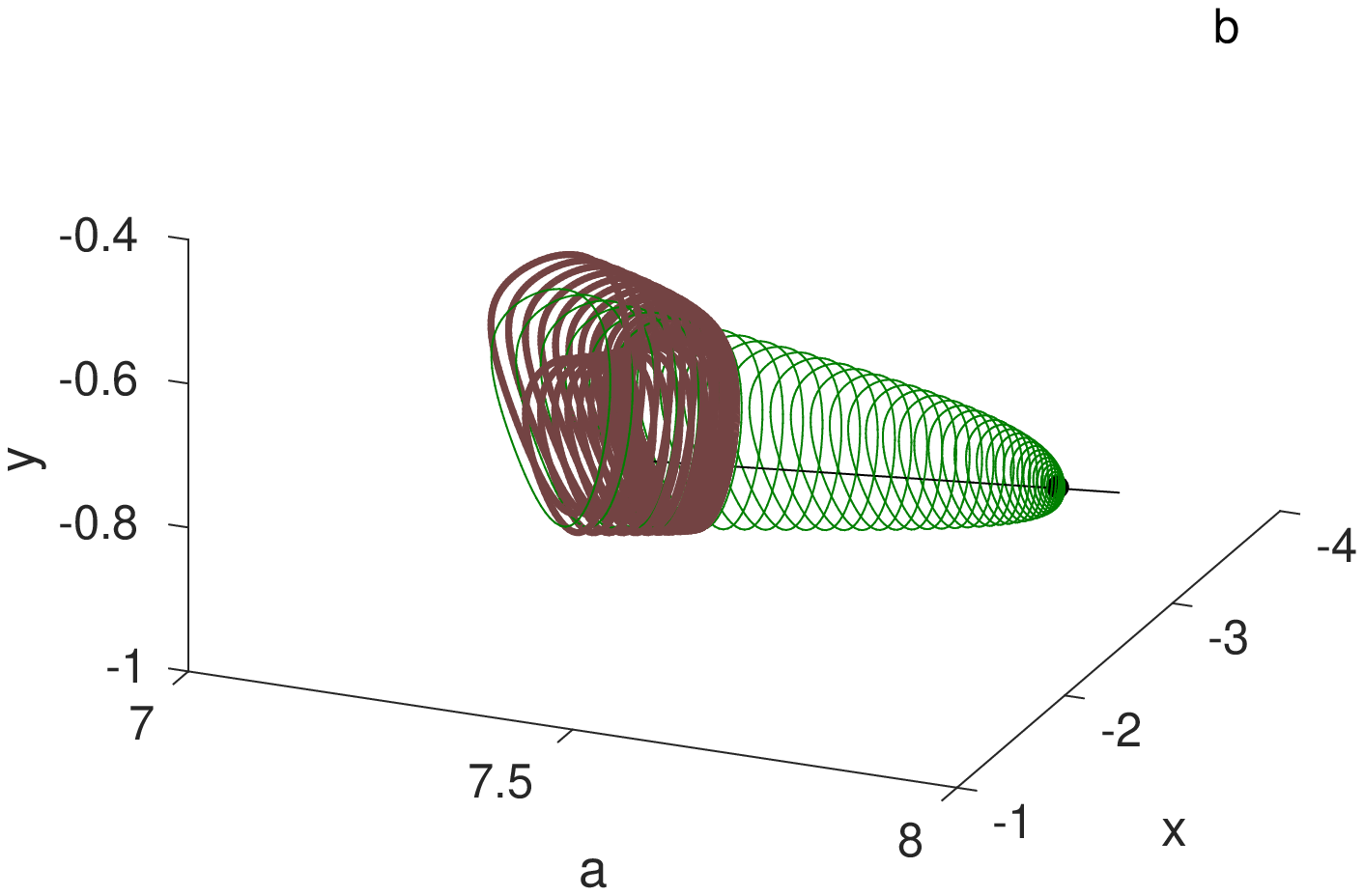}
 \caption{ (a) (Blue) Hopf bifurcation and (red) period doubling curves in the $(\delta,a)$ plane. The point near $(7.3945, 0.2077)$ labelled GH separates subcritical and supercritical bifurcations~\cite{matcont} along the Hopf curve. The black asterisk shows the $(\delta,a)$ parameter values $(7.3939,0.225423)$ studied extensively in this paper. These values are  below the period doubling curve. The parameters $[\rho,c,k] = [0.3224,2.3952,0.4032]$
(b) A bifurcation diagram with varying parameter $a$ showing a curve of equilibria (black) that passes through a Hopf bifurcation point (black dot) near $\delta = 7.9224$, the (green) family of periodic orbits emerging from the Hopf point and the  (brown) family of ``doubled'' periodic orbits emerging from a period doubling bifurcation near $\delta = 7.4438$. The parameters $[\delta,\rho,c,k] = [0.225423,0.3224,2.3952,0.4032]$}
\label{enso_e_bif}
\end{figure}

Several different types of attractors are found in the JT model. We are especially interested in those that have strong El Ni\~no events. Since the variable  $x$ represents  the difference in sea surface temperatures between the eastern and western tropical Pacific, these events are marked by small values of $|x|$. We used bifurcation analysis in our search for different types of attractors. 
The program MATCONT~ \cite{matcont} was used to locate points of Hopf 
bifurcation~\cite{GH83} and then to identify points of ``generalized Hopf'' bifurcations 
where  the bifurcation is neither subcritical nor supercritical. Starting at 
parameters with a supercritical Hopf bifurcation and varying parameter $\delta$, periodic 
orbits were continued until they lost stability at a period doubling bifurcation and then went through a period doubling cascade to a chaotic attractor~\cite{eckmann81}. Figure~\ref{enso_e_bif} displays results of these MATCONT calculations. 
Further increases in $\delta$ led to a regime with an 
MMO attractor. Decreasing $\delta$ showed that the ranges of $\delta$ with MMO 
and chaotic attractors overlapped. 
At the end of this section, we discuss the bifurcations that occur on the boundaries 
of the overlap region. 

\begin{figure}
\includegraphics[width=0.45\textwidth]{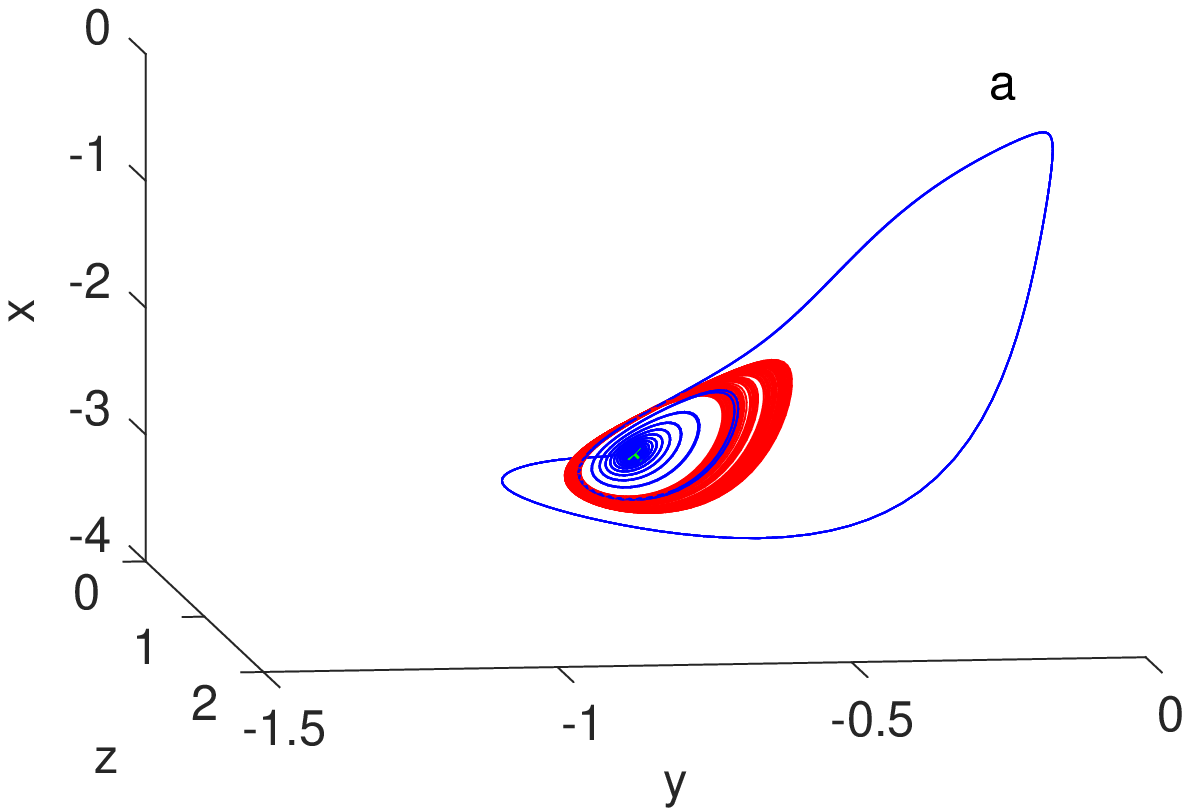}
\hspace{0.05\textwidth}
 \includegraphics[width=0.45\textwidth]{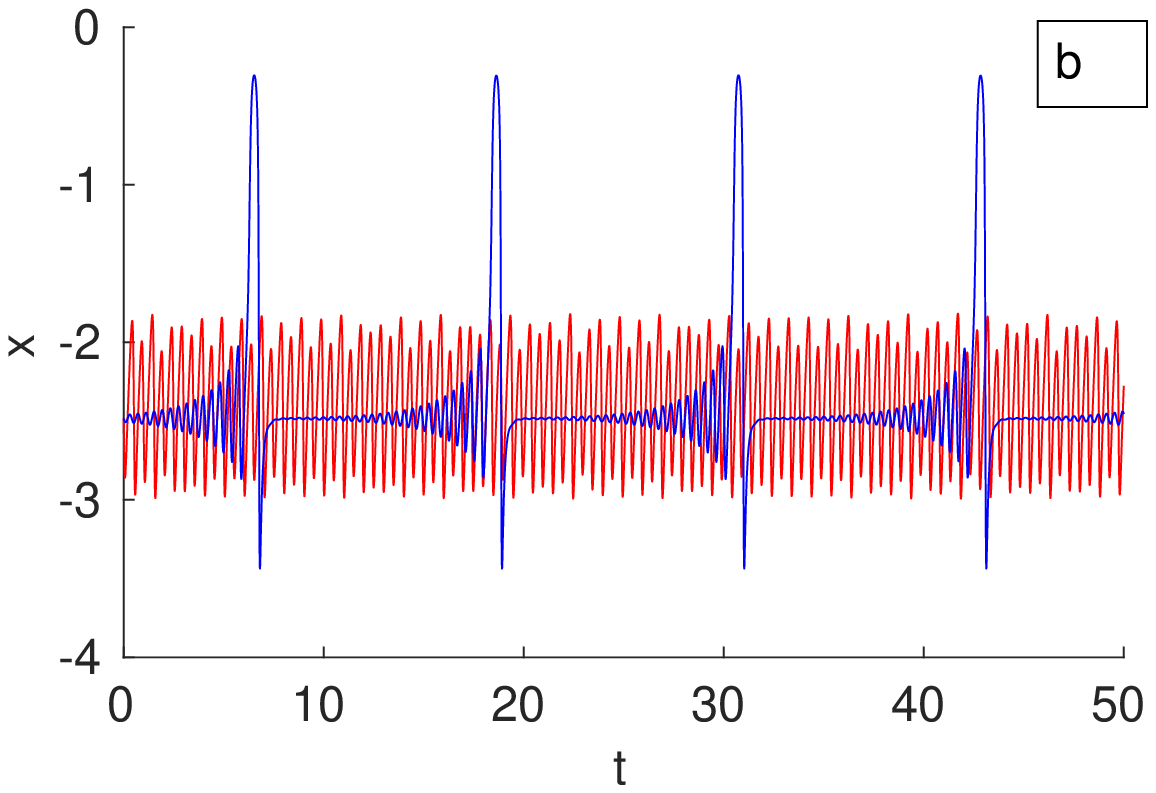}
  \caption{ (a) Phase portraits of a chaotic attractor (red) and an MMO periodic orbit (blue) that coexist at 
parameters  $[\delta,\rho,c,k,a] = [0.225423,0.3224,2.3952,0.4032,7.3939]$. 
The equilibrium point is a partially obscured blue dot. (b) Time series of the $x$ coordinate of the two trajectories.}
\label{enso_e_bistable_pp}
\end{figure}

The parameters
$$[\delta,\rho,c,k,a] = [0.225423,0.3224,2.3952,0.4032,7.3939]$$ 
in the overlap region are used here to study bistability and mode switching in the JT model.
Figure~\ref{enso_e_bistable_pp} shows two  trajectories of the model with these 
parameter values. The blue 
trajectory appears to be  a periodic mixed mode oscillation (MMO) with period 
approximately 12.1 years and a single intense El  Ni\~no event each period that follows 
a series of growing fast (sub-annual) oscillations. We show below that this trajectory is more complicated with small variations from one cycle to the next. The time scales in these dynamics are slightly shorter than those seen in observations of ENSO, but our purpose here is to describe a new fundamental mechanism that could potentially limit the predictability of strong El  Ni\~no events rather than finding an optimal quantitative fit of this highly reduced model to data.  The red 
trajectory is chaotic with oscillations whose amplitude variation is smaller 
than those of the MMO.

A central feature of the 
MMO cycle is an equilibrium point that is a 
saddle focus: the MMO approaches the 
equilibrium along its one dimensional stable manifold and then spirals away,
following its two dimensional unstable manifold. The number 
of oscillations  within the MMO cycle and their minimum amplitude depends upon how close it 
approaches the equilibrium. The growing amplitude oscillations of the MMO cycle 
terminate in a strong El Ni\~no event during which the sea surface temperatures of eastern 
and western tropical Pacific become almost the same. This El Ni\~no is followed rapidly by 
a La Ni\~na event in which the system ``recharges'' \cite{Jin1997},
reestablishing  higher sea surface temperatures and thermocline depth in the 
western than eastern Pacific. 
The  recharged system flows back toward the 
equilibrium where the cycle repeats.  

\begin{figure}
\includegraphics[width=0.45\textwidth]{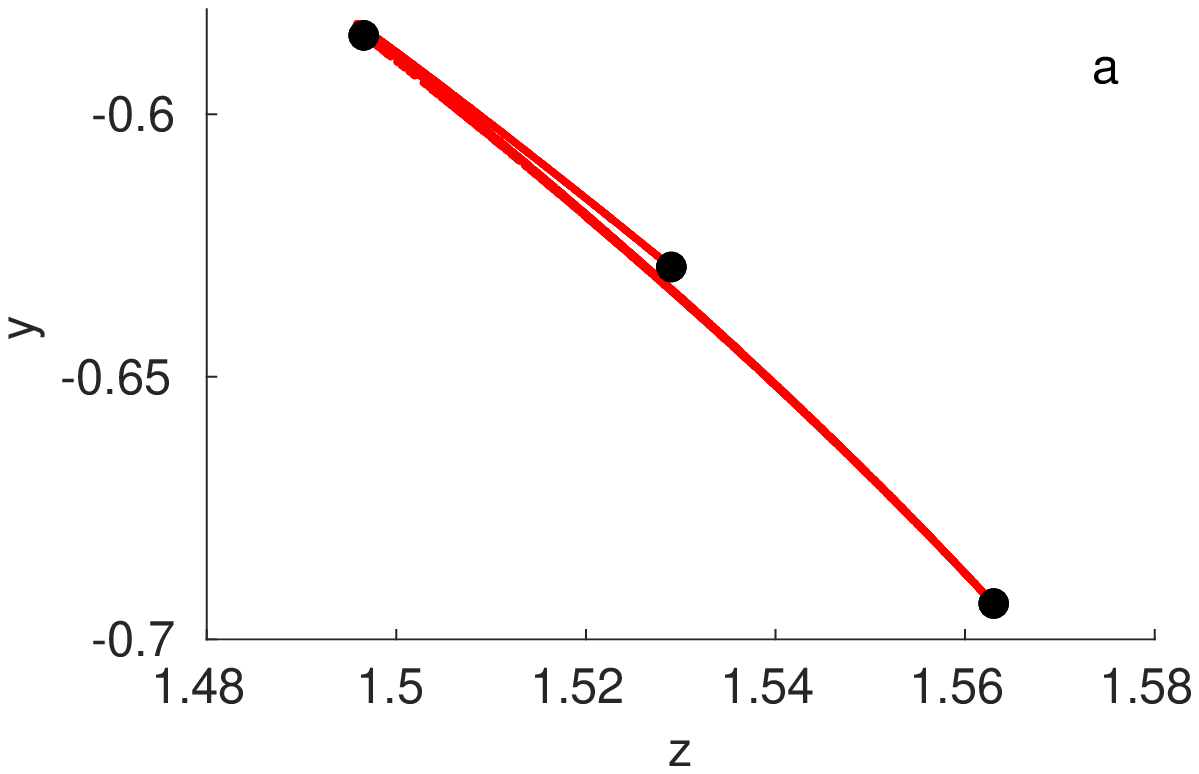}
\hspace{0.05\textwidth}
 \includegraphics[width=0.45\textwidth]{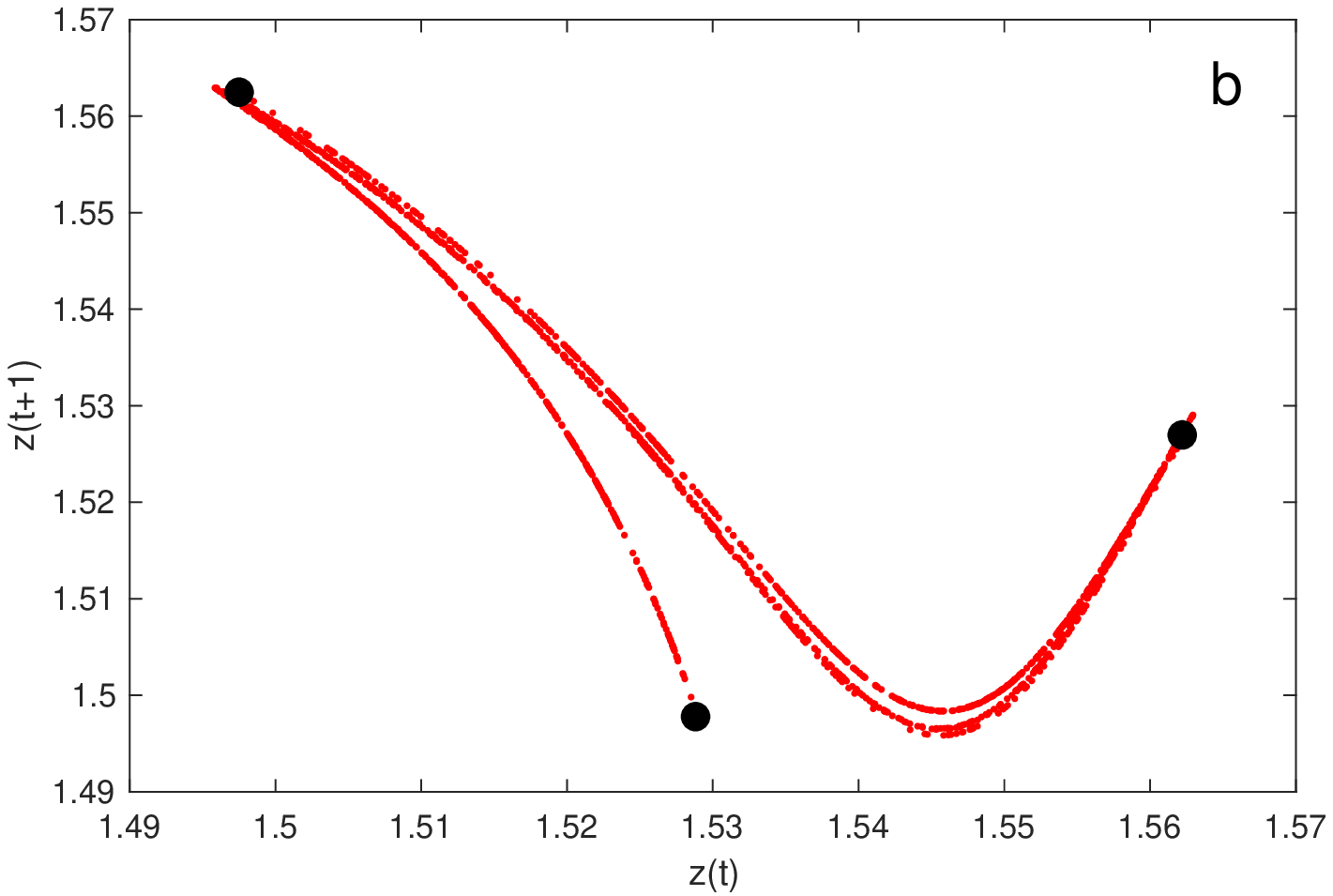}
  \caption{ (a) Intersections $(y_n,z_n)$ of a trajectory in  the chaotic attractor with the plane $x=x_{eq} \approx -2.4839$. Black dots show intersections of a  periodic orbit lying on the basin boundary with the cross-section. (b) $z_{n+1}$ against  $z_{n}$ for this trajectory.}
\label{enso_e_small}
\end{figure}

The chaotic attractor resembles those found in many other three dimensional 
vector fields. In contrast to strong El Ni\~no events, $x$ remains far from $0$.
As illustrated in Figure~\ref{enso_e_small}, 
this attractor can be analyzed by introducing a cross-section and studying 
its return map. Figure~\ref{enso_e_small}a shows the intersection of a trajectory 
with the cross-section $x=x_{eq}$, where $x_{eq} \approx -2.4839$ is the value of $x$ at the 
equilibrium point. Black filled circles in the figure are intersections of a periodic orbit with the cross-section. This periodic orbit appears to be on the boundary of the basin of attraction of the attractor, and it plays a central role in bifurcations of the attractor. The return map contracts in one 
direction, and stretches and folds in a second direction. Figure~\ref{enso_e_small}b
shows the folding by plotting the value
of $z$ at each return to the cross-section against the value of $z$ at the previous return. 
There is sensitivity to initial conditions within the chaotic attractor, but it is 
relatively mild: nearby initial conditions separate along the attractor without 
leaving it.

\begin{figure}
\centering
\includegraphics[width=0.9\textwidth]{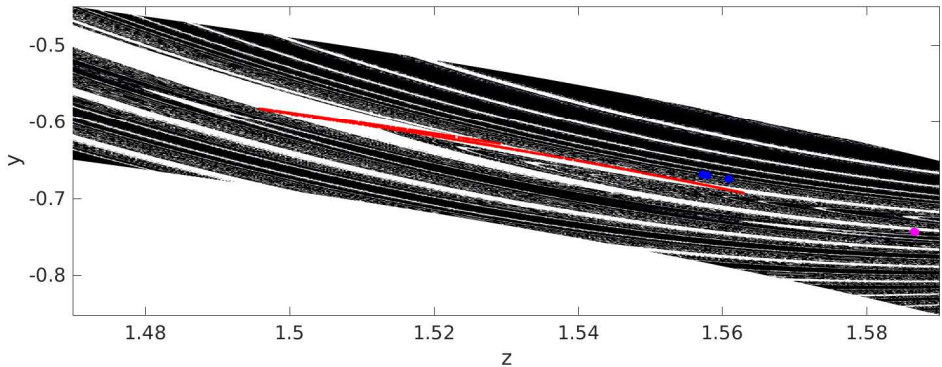}
  \caption{ The colored region is a 
500x500 grid that straddles the 
intersection of the chaotic attractor with the cross-section $x = x_{eq} \approx -2.8439$. Points 
are colored by whether their Trajectories of light colored points approach 
the chaotic attractor while trajectories of dark colored points approach 
the MMO attractor. The red markers and blue 
dots mark intersections of the chaotic attractor and MMO with the cross-section, 
respectively. The (overlapping) magenta markers locate initial points for the trajectories displayed in 
Figure~\ref{enso_e_sdic}a. 
}
\label{enso_e_imap}
\end{figure}

The basins of attraction of the two attractors are intertwined with each other 
in an intricate way, as illustrated in Figure~\ref{enso_e_imap}.  Trajectories starting 
in a $500 \times 500$ initial condition grid in the cross-section $x=x_{eq} \approx -2.4839$
were computed. 
The figure colors those initial  conditions that approach the chaotic attractor yellow and 
those that approach the MMO attractor blue. Intersections of the chaotic attractor with the 
cross-section are plotted red and the intersections of the MMO attractor with the cross-section 
are plotted blue. The attractors appear to have fractal basin boundaries~\cite{grebogi87}. 
Further evidence for this assertion was obtained by computing a saddle  periodic orbit 
$\Gamma$ with negative Floquet multipliers in the basin boundary. Were the basin boundary smooth, $\Gamma$ would separate its unstable manifold into two  pieces, one consisting of trajectories that approach the MMO attractor and one consisting of trajectories that approach the chaotic attractor. However, a trajectory with initial condition on one side of the unstable manifold returns on the opposite side because the multipliers are negative. Consequently, the local unstable manifold of $\Gamma$ is a one sided, non-orientable M\"obius strip that cannot constitute the common boundary of two different attractor basins. Instead, the two basins of attraction interleaved thin layers and each approaches $\Gamma$ from both sides of its stable manifold. Thus, 
there are large sets in which it is difficult to predict which basin a 
chosen initial condition will belong to.
Both attractors are very close to their basin boundaries, so carefully 
chosen perturbations of very small magnitude can induce trajectories to flow from 
the chaotic to the MMO attractor and vice versa. To explore mode switching in this model, 
we investigated regions where nearby trajectories separate quickly.

\begin{figure}
 \includegraphics[width=0.45\textwidth]{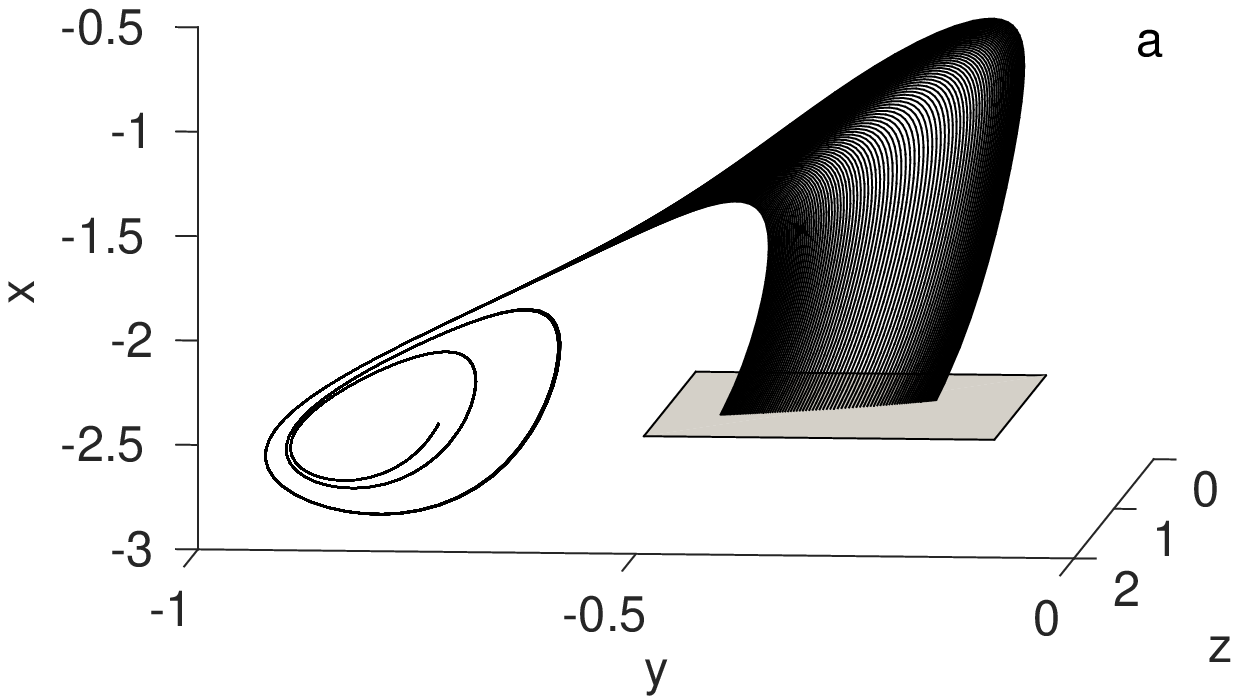}
 \hspace{0.05\textwidth}
 \includegraphics[width=0.45\textwidth]{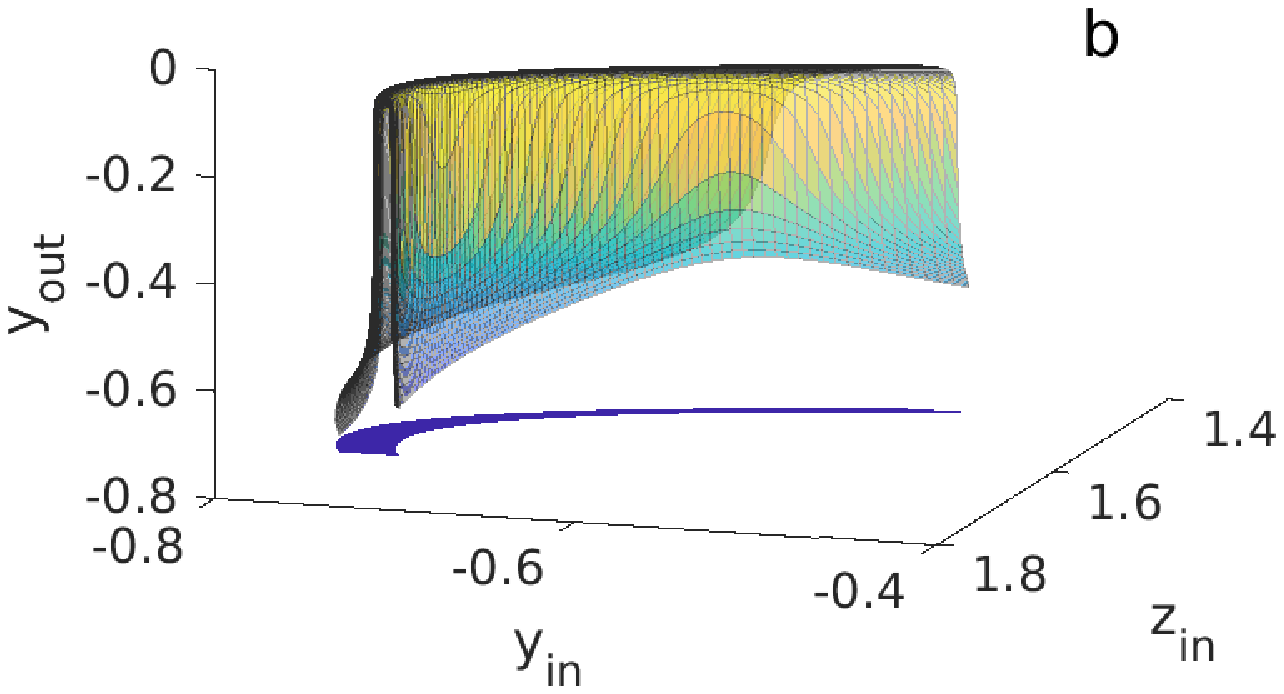}
   
\caption{(a): Ensemble of  100 trajectory segments in the cross 
section $x = x_{eq} \approx -2.4839$ are plotted to their third return to this cross-section, drawn gray.
Initial conditions lie along a curve of length approximately $0.0005$. The 
final points of the segments lie on a curve of length larger than $0.4$.
(b): The image of the $y$-coordinate under the third iterate of the 
return map  of a 
thin strip of initial conditions (dark blue) of width 0.035 in the cross 
section $x= x_{eq}$. Points in the middle of the strip make a large amplitude 
excursion representing a strong El Ni\~no. The map stretches the $y$-coordinate 
by a factor of order $10^3$ where the returns have intermediate amplitude between those on the ridge and the strip boundaries.
}
\label{enso_e_sdic} 
  \end{figure}

\begin{figure}
\centering
\includegraphics[width=0.9\textwidth]{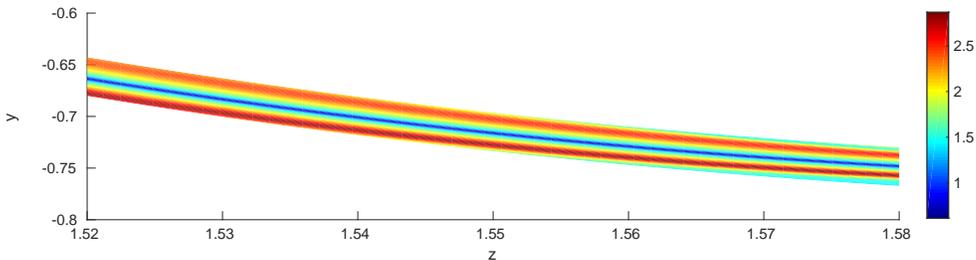}

  \caption{ The largest finite time Lyapunov exponent for trajectory segments of 
length five. Initial conditions are a $50 \times 50$ grid in the plane $x= -2.4839$ and 
values of the FTLE are plotted on a base 10 log scale. Points in the middle of the strip 
make a large amplitude excursion representing a strong El Ni\~no event. }

\label{enso_e_ftle}
\end{figure}

Figure \ref{enso_e_sdic} visualizes the rapid separation of trajectories. Figure \ref{enso_e_sdic}a shows a collection of 100 trajectory segments with initial conditions in the cross-section $x = x_{eq}$. These initial conditions were chosen along a curve of length approximately $0.0005$ whose image under the third iterate of the return map appears to stretch between the two attractors. Figure \ref{enso_e_sdic}b shows the $y$ coordinate for the third return of points in a region  $W$ similar to  the region displayed in Figure~\ref{enso_e_imap}. In the middle of this strip, there is a ``ridge'' of points for which the $y$ coordinate of the return is close to $-0.1$. These are points whose third return is associated with a strong El Ni\~no event. The sides of the ridge are points where the third iterate of the return map has large stretching. This can be quantified by computing  the finite time Lyapunov 
exponent or FTLE \cite{Haller2001} of trajectories in $W$. The FTLE is the largest singular value of the variational equations of the flow map along a trajectory segment of specified length. It measures the maximal stretching of infinitesimally close trajectories along the trajectory segment. Figure  \ref{enso_e_ftle} plots this quantity on a log scale in a strip of initial conditions in the plane $x= -2.4839$ and trajectory segments of time 
length $ \log_{10} 0.7 \approx 5$.  Values are on a base 10 log scale. The center of the strip, indicated  by its blue color,  consists of 
points in the MMO basin of attraction and their FTLEs are approximately $5$. Trajectory
segments along two flanking strips have FTLE two orders of magnitude larger,
demonstrating the high sensitivity to initial conditions there.

These properties
contrast with the one dimensional map $g$ discussed in the previous section. As
the parameter $a$ of that example increases to its critical value $3\sqrt{3}/2$, 
the basins of 
two attractors touch at a single point and the two attractors merge into one. In the JT model,
bifurcations of  the MMO and chaotic attractors occur at different parameter values and the two 
attractors do not merge to become one. 
Apparently, when one of the invariant sets ceases to be an attractor, most points that were in its basin of attraction have trajectories that now flow to  the other attractor. We investigated how this happens as the parameter $a$ is varied.  

In the case of 
the chaotic attractor of the JT model, a \emph{boundary crisis bifurcation}~\cite{GOY_crisis}  occurs where the attractor
becomes a non-attracting invariant set. As $a$ decreases to a value close to
$7.3938587$ the periodic orbit on the basin boundary displayed in Figure~\ref{enso_e_small} approaches the attractor. For smaller values of $a$, a chaotic invariant set persists, but it is no longer attracting. Trajectories can ``escape'' from the former attractor basin along the unstable manifold of the periodic orbit. 

The bifurcations associated with the MMO attractors also involve a boundary crisis. The MMO attractor for $a=7.3939$ comes close to the equilibrium point which is a saddle focus with 
eigenvalues approximately $-1.46$ and $0.127 \pm 4.47$. Because the ratio of the negative eigenvalue to the real part of the complex eigenvalues has large magnitude, a large volume of phase space is drawn toward the stable manifold of the equilibrium before their trajectories  spiral away slowly along the two dimensional unstable manifold of the equilibrium. Moreover,  trajectories starting close to the equilibrium in the unstable manifold increase their distance from it by a factor of only about $1.2$ per revolution as determined by the ratio of the real and imaginary parts of the complex eigenvalues.  Nonetheless, motion along the unstable manifold has the potential to disperse trajectories through a large part of the phase space. Some of the trajectories in the unstable manifold undergo strong El Ni\~no  events and approach the MMO attractor. 

\begin{figure}
 \includegraphics[width=0.45\textwidth]{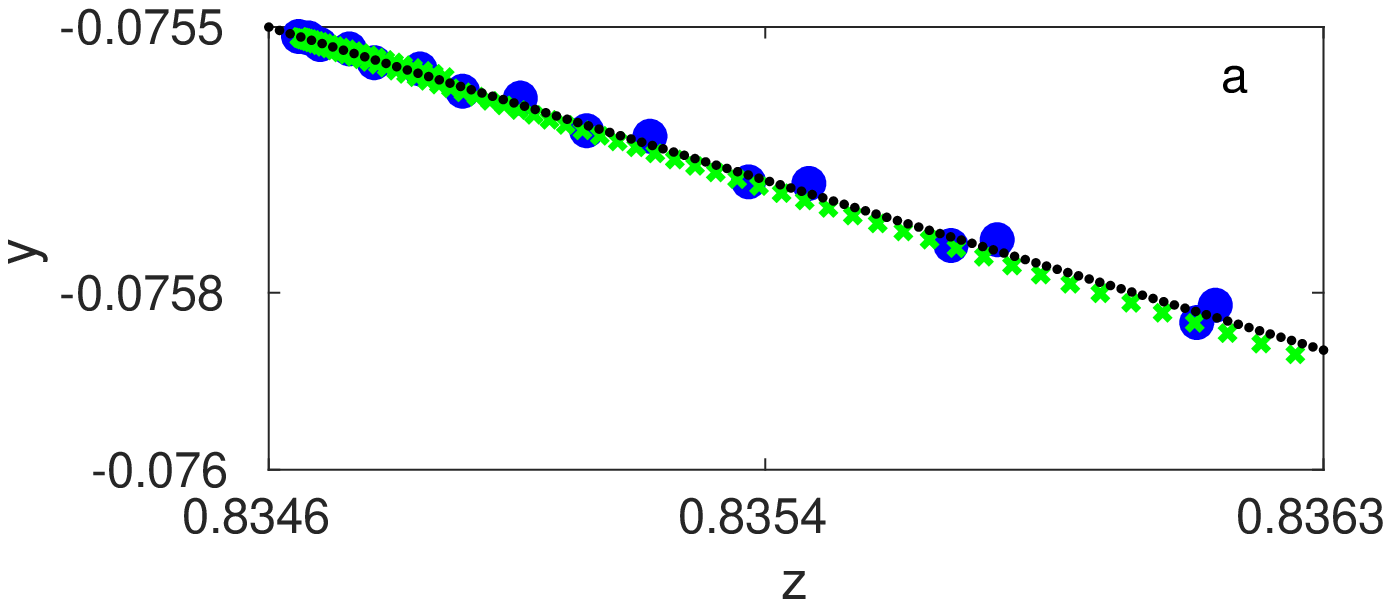}
 \hspace{0.05\textwidth}
 \includegraphics[width=0.45\textwidth]{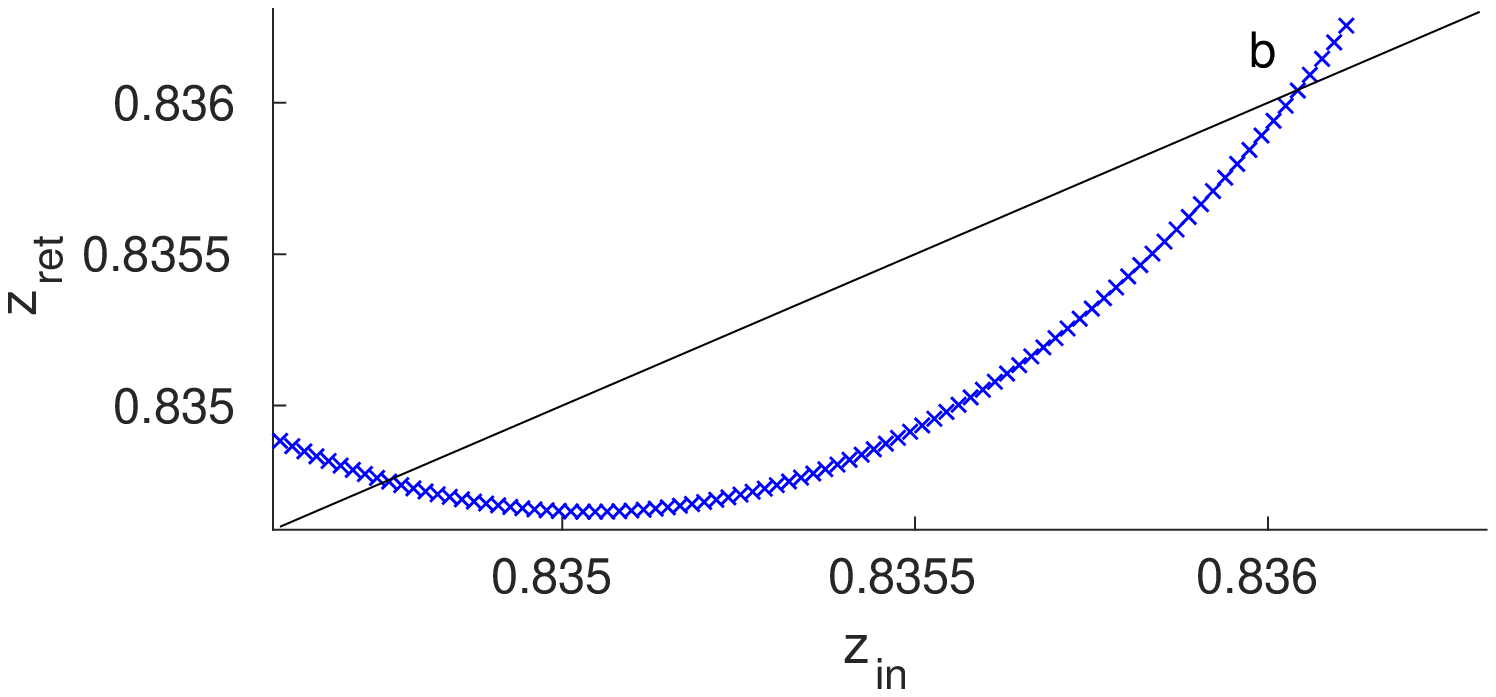}
 \\
  \includegraphics[width=0.45\textwidth]{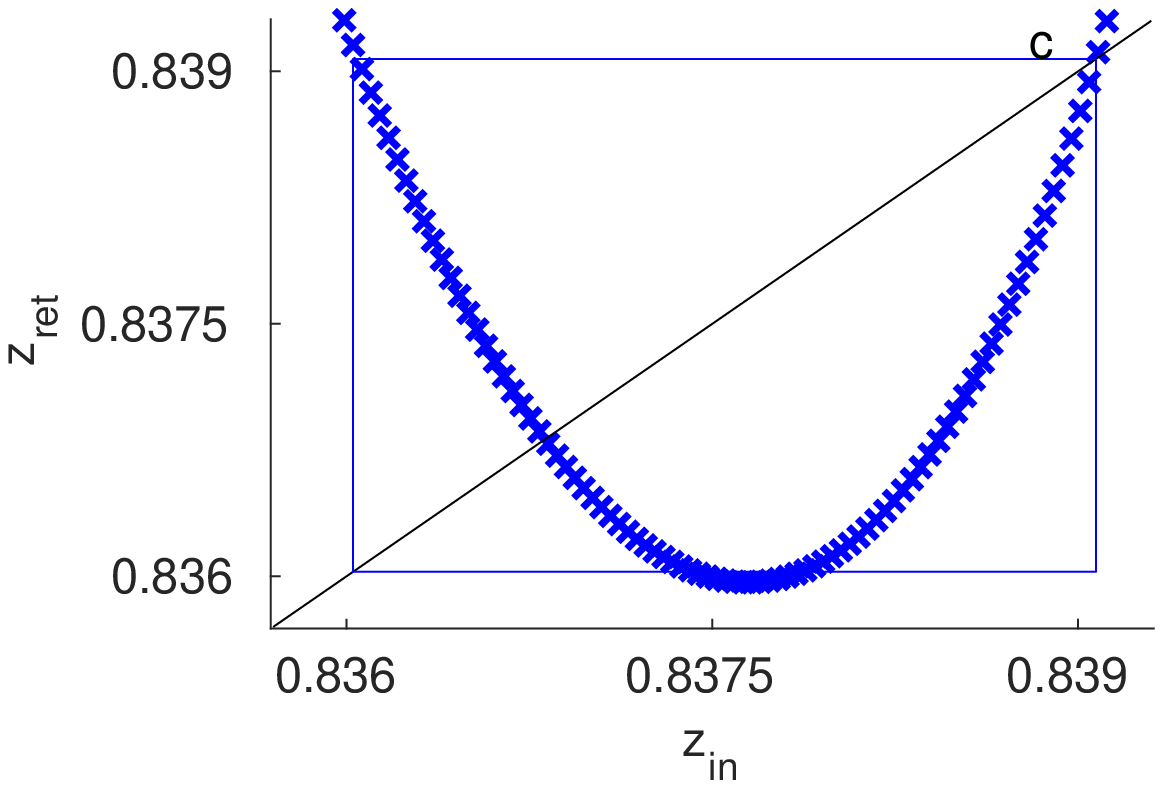}
 \hspace{0.05\textwidth}
 \includegraphics[width=0.45\textwidth]{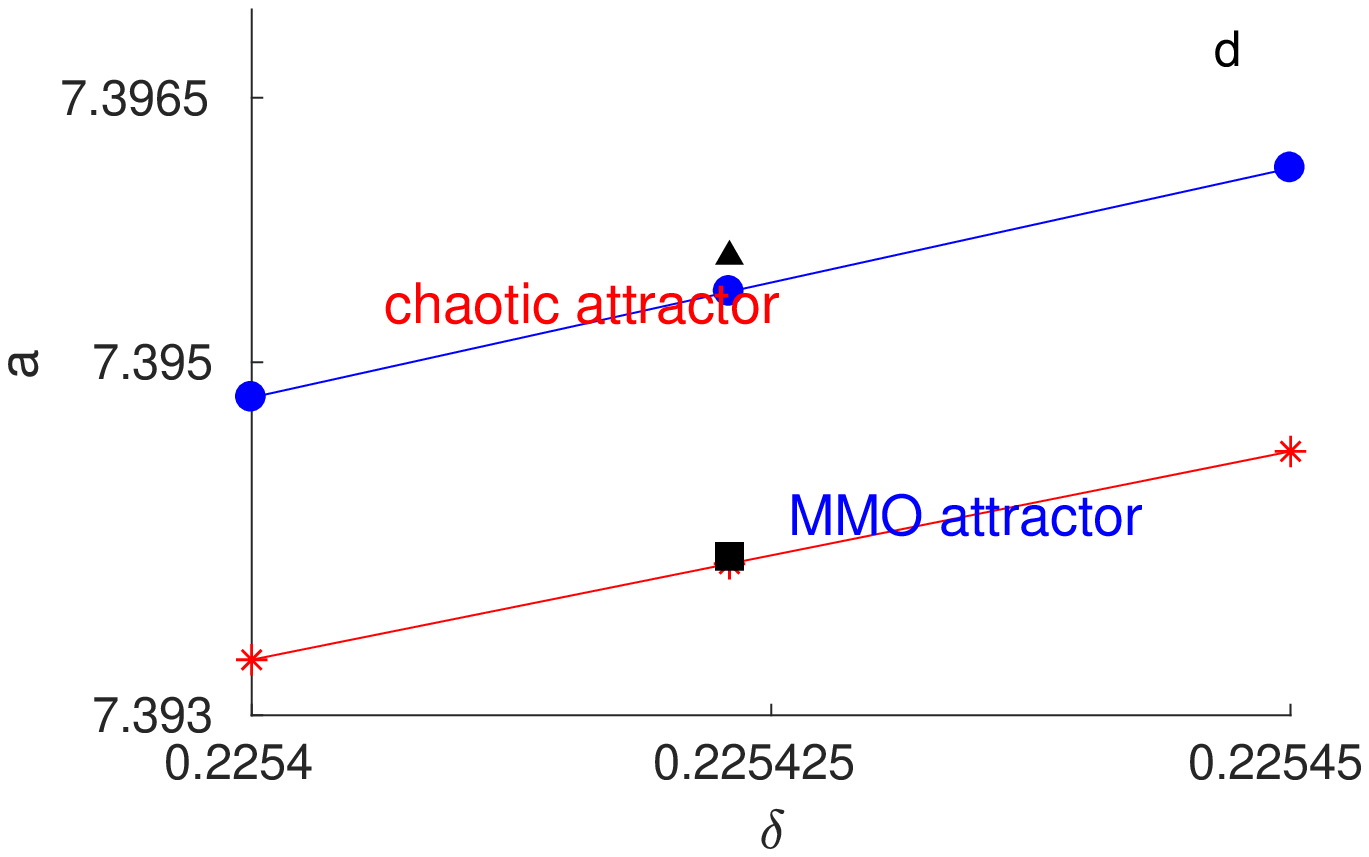}

\caption{(a) Intersections of the unstable manifold of the equilibrium with the cross-section $x=-1.5$ are plotted as blue filled circles. A quadratic fit to these points is plotted as black dots. Images under the return map for the black points are plotted as green x's.
(b) Return map values of the $z$-coordinate plotted against their initial values on the black curve of subplot (a). The black line is $z_{ret} = z_{in}$. The parameters 
$(\delta,a) = (0.225423,7.393)$.
The interval $[0.835,0.837]$ is mapped into itself in a unimodal  manner.
(c) The approximate one dimensional return map for $(\delta,a) = (0.225423,7.3956)$. Note that this return map has a non-attracting horseshoe.
(d) Approximate bifurcation boundaries for the two types of attractors. Chaotic attractors are present in the region that lies above the red curve marked by asterisks. Along this curve, a periodic orbit touches the attractor in a boundary crisis. MMO attractors are present in the region below the blue curve marked by filled circles. Crossing this curve from below to above, MMO attractors become invariant sets that intersect the cross-section $x=1.5$ in a non-attracting horseshoe. The bistability region where both types of attractors exist is the strip between the two curves. The parameters $(\delta,a) = (0.225423,7.3939)$ inside the strip, but near the chaotic bifurcation boundary, are displayed as a black square. Figure 3 displays points of the periodic orbit that collides with the chaotic attractor in a boundary crisis. The parameters $(\delta,a) = (0.225423,7.3956)$ outside the bistability strip are marked by a black triangle. Subplot (c) displays the non-attracting horseshoe for the return map to $x=1.5$ for these parameters.
}
\label{enso_e_ret_73939} 
  \end{figure}

We examine next a return map to the cross-section $x=-1.5$, which we set (arbitrarily) as a threshold for strong El Ni\~no events. This value is larger than the maximum of $x$ in the chaotic attractor (see Figure \ref{enso_e_bistable_pp}.)  Figure~\ref{enso_e_ret_73939}a plots intersections of the unstable manifold of the equilibrium with this cross-section as large blue dots. There is a fold that occurs near $z=0.835$ dividing the intersection into two branches that are close together. To investigate the next return of the unstable manifold to the cross-section, we approximate the intersection by a quadratic curve, plotted as black dots in the figure. The returns of the black points are plotted as green x's. This image lies close to the intersection of the unstable manifold of the equilibrium with the cross-section, so we regard the return map as approximately one dimensional. Figure~\ref{enso_e_ret_73939}b shows the graph of this one dimensional map, parametrized by the z coordinates of domain and range. There is an interval (roughly $z \in [0.835,0.837]$) mapped into itself. A further computation using Newton's method on the return map locates a fixed point  $p$ at $(y,z) \approx (-0.07566  ,0.83542)$ with one eigenvalue approximately -1.245 and the other smaller than $10^{-6}$, confirming that the return map is approximately rank 1. Since the magnitude of one eigenvalue is larger than $1$, the periodic orbit containing $p$ is not stable, precluding that it is the MMO attractor. Nonetheless, the MMO attractor lies close to this periodic orbit. As the parameter $a$ varies, the one dimensional approximations of the return map to the cross-section $z=1.5$ undergo a ``full'' set of bifurcations for unimodal maps. The first bifurcation of the return map is a saddle-node bifurcation of a fixed point in the $a$ interval $[7.3915,7.3939]$. This bifurcation produces a stable periodic MMO with a single strong El Ni\~no  event each cycle. This is followed by a cascade of period doubling bifurcations and chaotic MMOs as $a$ increases~\cite{eckmann81}. At  $a=7.3956$, the return map has a horseshoe; i.e, a non-attracting invariant set topologically equivalent to a shift on two symbols~\cite{GH83}. This is illustrated in Figure~\ref{enso_e_ret_73939}c. All of these MMO attractors have only small modulations of an MMO periodic orbit and make repeated passages near the equilibrium point, producing strong El Ni\~no  events of similar frequency.

Having identified bifurcations on the boundaries of the strip of the $(\delta,a)$ parameter space where there is bistability. Values of $\delta$ were selected and then a divide and conquer strategy was utilized with varying $a$. In the region with chaotic attractors, the distance between the attractor chaotic attractors and the periodic orbit shown in Figure 3 indicates the distance to the bifurcation boundary.  The scale of Figure 3 does not allow one to see this distance, but it is clearly visible for values of $a$ farther from the boundary. Newton's method was applied to the return map in order to compute the periodic orbits with high precision. The extent of the attractor was determined by computing intersections of along trajectory segment with the cross-section. On crossing the bifurcation curve, we observed a discontinuous increase in the set of intersections  of trajectory and cross-section. For the MMO attractors, we computed approximate return maps like those displayed in Figure~\ref{enso_e_ret_73939}b,c. The boundary crisis occurs where the image of the minimum value of the return map is its fixed point of positive slope. This is similar to the bifurcation displayed for the one dimensional map $g$ in Section 2.  In Figure~\ref{enso_e_ret_73939}c, the image of the minimum value is clearly larger than the fixed point, so that the map has a horseshoe. Figure~\ref{enso_e_ret_73939}d displays three points on each of the two bifurcation boundaries.This subplot also shows the parameter values $(\delta,a) = (0.225423,7.3939)$ that we have studied intensively as well as the parameters $(\delta,a) = (0.225423,7.3956)$ that lie just outside the bistability strip. These (minimal) computations show that the strip of bistability is thin but not exceedingly so. On the scale of the figure, the boundaries are close to parallel, suggesting that the strip may be quite long.

\begin{figure}
 \includegraphics[width=0.45\textwidth]{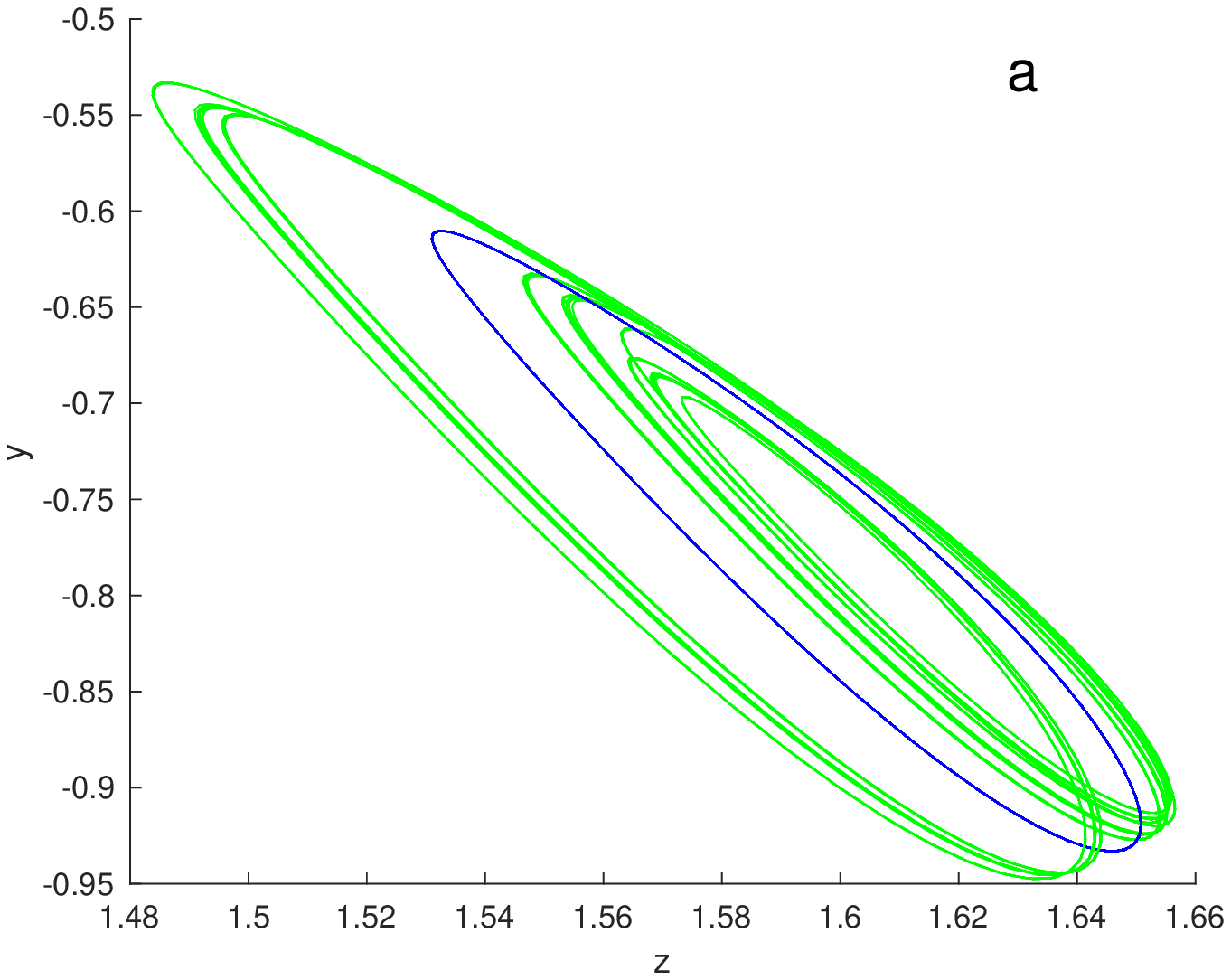}
 \hspace{0.05\textwidth}
 \includegraphics[width=0.45\textwidth]{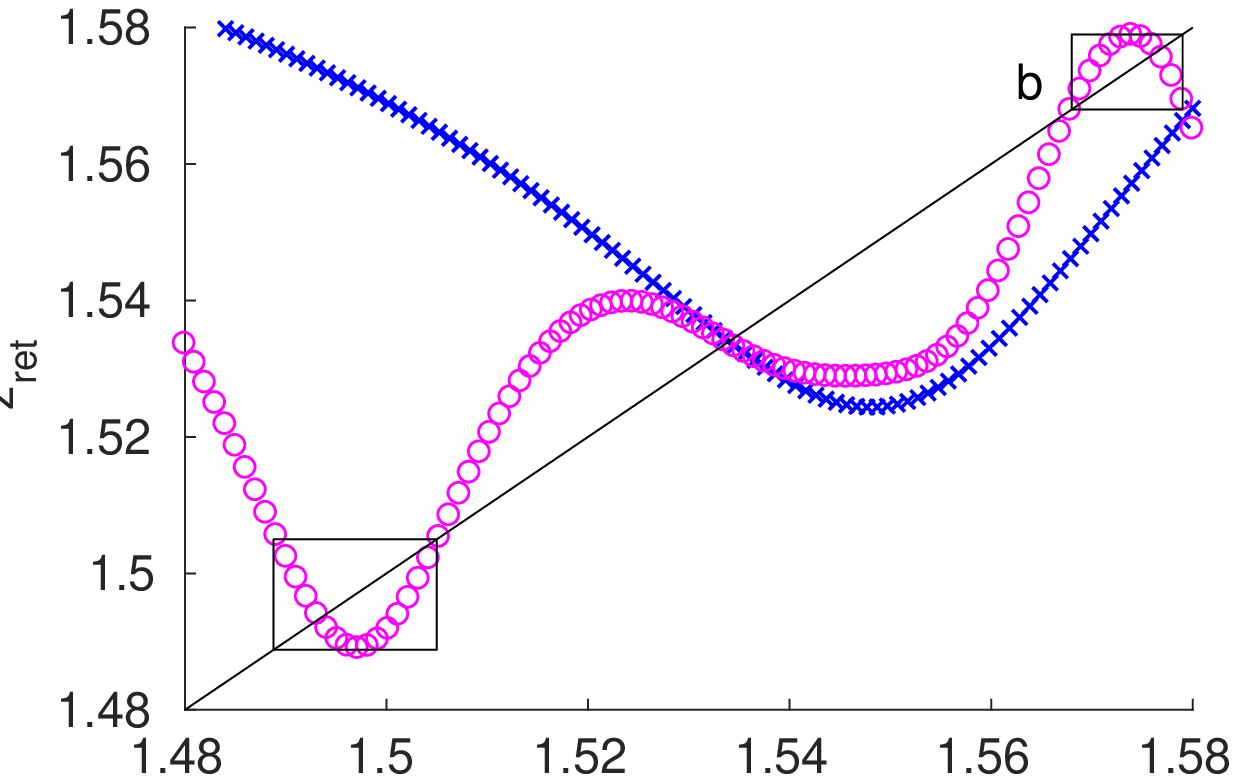}
   
\caption{(a): Two trajectories for parameter value $a=7.453$, a periodic orbit (blue) and a chaotic attractor (green).
(b): Images of the $z$-coordinate under 
return maps  for a curve of initial conditions that lie close to the image of the return map. Blue x's are images of the return map; magenta circles are images of the third iterate. The third iterate has subintervals (shown by black boxes) that are mapped into themselves in addition to a stable fixed point near 1.534.
}
\label{enso_e_bistable_7453} 
  \end{figure}

These explorations of the parameter space of the JT model only begin to probe the dynamics found in the model. There are also trajectories in the unstable manifold of the equilibrium that do not make strong El Ni\~no  excursions. The boundary between the regions that do and do not make these excursions occurs in the strips of points with large FTLE described above. We expect to find chaotic, non-attracting invariant sets in these regions. For values of  $a>7.453$, the unstable manifold of the equilibrium no longer seems to contain trajectories that have strong El Ni\~no events, precluding the existence of MMO attractors like those we have described. This is also a situation that bears resemblance to phenomena observed in the Koper model~\cite{koper1992}. There, MMO trajectories are observed when the unstable manifold of the equilibrium intersects a repelling slow manifold~\cite{mmoreview} of a two-time scale system with two slow variables and one fast variable. Tangency of these two manifolds locates the parameter space boundary  between regions with and without MMOs. Similarly here, we observe strong El Ni\~no  events when the unstable manifold of the equilibrium intersects the strip with large FTLE. Even when there are no MMO attractors, bistability may be present. Figure~\ref{enso_e_bistable_7453} illustrates bistability between a periodic orbit and a chaotic attractor. These examples exhibit diverse types of  ENSO dynamics in the JT model, but we do pursue them further in this paper. Our main goal is to demonstrate unpredictability where mode switching between attractors with and without strong El Ni\~no events occurs abruptly. This is the subject of the next section.

\section*{Mode Switching}

\begin{figure}
\centering
  \includegraphics[height=0.25\textwidth]{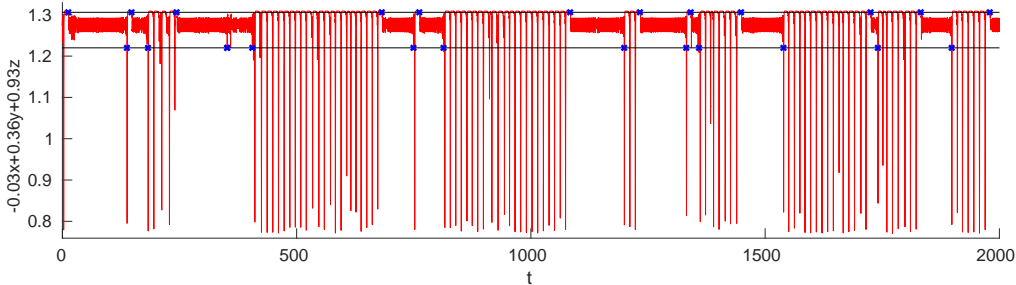}
  \caption{Time series of the function $g(x,y,z) \approx -0.0303x+0.3600y+0.9325z$ on a 
  portion of a trajectory computed with $a = 7.3939+0.002\sin(1.8 t)$ and $[\delta,\rho,c,k] = [0.225423,0.3224,2.3952,0.4032]$.  The units of time are years. The crosses mark a subset of points where the trajectory crosses the thresholds $g = 1.306$ and $g = 1.22$ as described in the text. The oscillations of $g$ lie between these thresholds during chaotic epochs. 
}
\label{enso_f_ts}
\end{figure}

The model of Jin and Timmermann \cite{tim03} includes terms representing annual
forcing from the seasonal cycle of solar insolation and stochastic forcing (e.g., due to rapidly varying wind stress on 
the sea surface). We found that the addition of either of these components to our 
simulations is capable of inducing mode switching between the two attractors discussed
in the previous section. Figure~\ref{enso_f_ts} illustrates a simulation of the periodically
forced system. The parameter $a$ was made a sinusoidally varying function of time
with period one year and amplitude $0.002$. Even at this small forcing amplitude, the system switches erratically between two almost invariant sets that approximate the attractors of the unforced model. The figure shows a representative time series of 2000 years selected from a simulation of 100,000 years. To visualize the almost invariant sets in a time series, we defined an ``observable'' $g$ to be the projection onto the line orthogonal to the eigenspace  of the unstable complex eigenvalues at the equilibrium point. The oscillations of $g$ near the chaotic attractor have much smaller magnitude than those of the MMO attractor. The epochs during which the trajectory is close to one of the two attractors are readily apparent in the time series of $g$. The switching between modes is quite abrupt, and the duration of the epochs is irregular. We sought to investigate the distribution and independence of these durations, but found it difficult to define precise phase space boundaries to use in segmenting the trajectory into epochs automatically. 

We set two thresholds for the function $g$ and computed the times when $g$ decreased through a threshold. The values of the thresholds were chosen to be intermediate between the range of oscillations of $g$ in the chaotic and MMO attractors. They are plotted as black lines in Figure~\ref{enso_f_ts}. Most of the chaotic attractor epochs lie between two successive events, with the upper threshold crossed at the beginning of the epoch and the lower threshold crossed at the end of the epoch. Moreover, the duration of most of these epochs is substantially larger than the roughly decadal duration between strong El Ni\~no events during their epochs. This prompted us to select pairs of successive events separated by a time duration longer than 15 years as bracketing chaotic epochs. In most cases, blue crosses on the upper threshold mark the onset of chaotic epochs, while those on the lower threshold mark the termination of the epochs. The complements of these time intervals were regarded as strong El Ni\~no epochs. In some instances, such as a time interval near $t=20200$, there are shorter oscillations in the chaotic range that have not been identified as chaotic epochs.

\begin{figure}
 \includegraphics[width=0.45\textwidth]{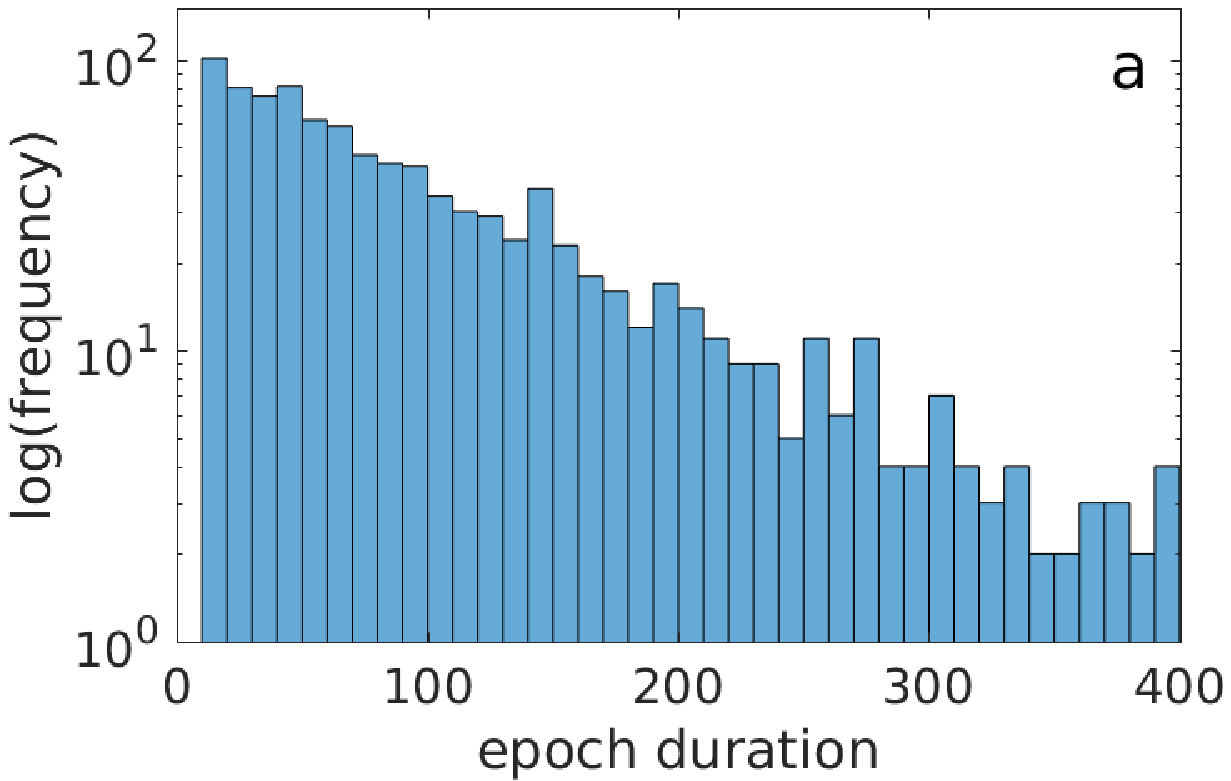}
 \hspace{0.05\textwidth}
 \includegraphics[width=0.45\textwidth]{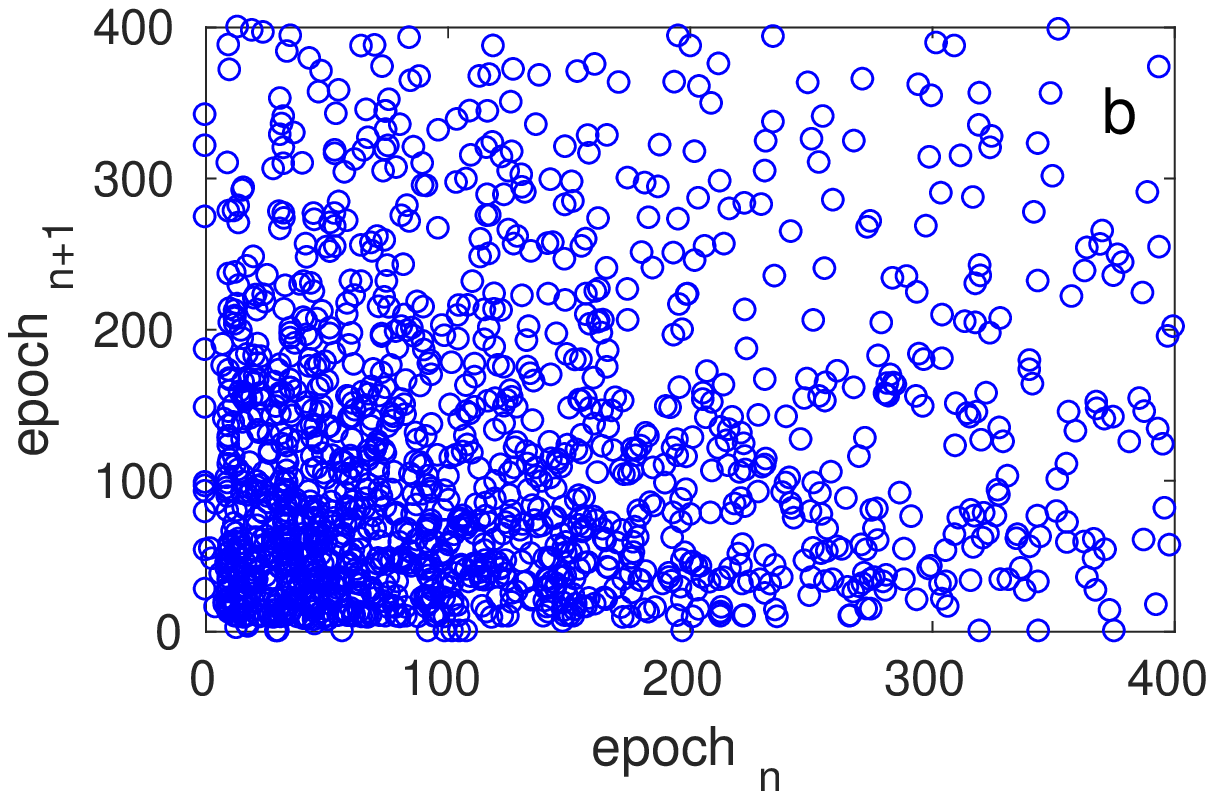}
   
\caption{(a): Histogram of epoch durations smaller than 400, sorted into bins of width 10 years. The number of occurrences in each bin are displayed on a log scale to highlight the resemblance to an exponential frequency distribution.
(b): Scatterplot of pairs of successive epoch durations smaller than 400. 
}
\label{enso_f_hist} 
  \end{figure}

Figure~\ref{enso_f_hist}(a) displays a histogram of the number of epochs of duration less than 400 years in bins of width 10. (The longest epoch has a duration of approximately 540 years.) These durations are quite long compared to those  observed in simulations of large climate models~\cite{Wittenberg2014}. However, we chose a very small amplitude for the periodic component of $a$ to make it easier to segment a trajectory into epochs. With larger magnitude of the periodic variation of $a$, the epochs are shorter and our criteria for segmenting a trajectory into epochs breaks down. The distribution of  epoch durations is roughly exponential. Figure~\ref{enso_f_hist}(b) is a scatter plot which shows pairs of successive epoch durations. Apart from higher density of shorter durations, there is no apparent pattern. We conclude that epoch durations and ENSO regime predictability in this JT model with small annual forcing are unpredictable~\cite{timmermann2006predictability}. 

\section*{Discussion}
Reconstructions of past climate from historical data indicate that ENSO is 
highly variable with long periods of larger and smaller variance on decadal 
time scales~\cite{cobb03,McGregorTimmermannClimateofPast2013}. Simulations 
of global climate models reach similar conclusions~\cite{Wittenberg2014}. 
There are many potential explanations for why ENSO is so variable. 
We have demonstrated here that interactions of only a few degrees of freedom in the 
highly reduced Jin-Timmermann model can produce unpredictability of strong 
El Ni\~no events and a new type of ENSO complexity on a decadal time scale. This is a stronger statement than 
saying that that ENSO is chaotic. It places a time scale on the sensitive dependence
on initial conditions. More specifically, our results show that epochs of relatively 
regular cycles of strong El  Ni\~no events like those observed during the past century
might not continue. If they do cease, they may resume at any time. 

We compare this unpredictability of ENSO regimes with ensemble weather forecasting. Ensemble forecasting is applied in situations where one expects sensitive dependence to initial conditions; i.e., nearby initial conditions yield trajectories that diverge from one another. This divergence happens gradually, so forecasting gives predictions with uncertainty that grows with time. Here, switching is abrupt and the time scale for predicting when a mode switch will take  place seems to be comparable to the decadal time scale of strong El  Ni\~no events. 
Rapid divergence of nearby  
trajectories from one another  is concentrated in small regions of the state  
space of the model. If we understand where this divergence takes place, then we can 
exploit that knowledge in our forecasts. We can also use what we have learned about
the complex dynamics of the low dimensional JT model as a guide to gain deeper insight into the coupled oceanic and atmospheric  dynamics of ENSO in large climate models.  Whether the large scale dynamics of these large models have low dimensional attractors and almost invariant sets displaying different ENSO dynamics remain unanswered questions. 

{\bf Acknowledgments:}  Henk Dijkstra's work was partially supported by a Mary Shepard B. Upson Visiting Professor position at the College of Engineering of  Cornell University, Ithaca, 
NY. He thanks for Prof. Paul Steen (Cornell University) for being his host and the many 
interesting discussions. Axel Timmermann was supported by the Institute for Basic Science (project code IBS-R028-D1).

\bibliography{enso_predictability_10_06_2017}
\bibliographystyle{plain}

\end{document}